%% file: isoperimetry_ps.tex
\documentclass[english,12pt]{smfart}
  \usepackage{macros_anglais}
  \usepackage{dsfont}

\title[Quantitative isoperimetry]{Sharp quantitative isoperimetric inequalities
in the L\textsuperscript{1} minkowski plane}
\author{Beno\^{\i}t Kloeckner}

\begin{document}

\maketitle

An isoperimetric inequality bounds from below the perimeter
of a domain in terms of its area. A \emph{quantitative}
isoperimetric inequality is a stability result: it bounds
from above the distance to an isoperimetric minimizer in terms of the isoperimetric
deficit. In other words, it measures how close to a minimizer an almost optimal set must be.

The euclidean quantitative isoperimetric inequality has been thoroughly studied, in particular in
\cite{Hall} and \cite{Fusco-Maggi-Pratelli}, but the $L^1$ case has 
drawn much less attention.

In this note we prove two quantitative isoperimetric inequality
in the $L^1$ Minkowski plane with sharp constants and determine the extremal domains
for one of them. It is usually (but not here) difficult to determine
the extremal domains in a quantitative isoperimetric inequality: the only
such kown result is for the Euclidean plane, due to Nitsch \cite{Nitsch}.

\section{Statement of the results}

We consider the plane $\mR^2$ endowed with the $L^1$ metric:
$$|(x_1,x_2)-(y_1,y_2)|=|x_1-y_1|+|x_2-y_2|.$$

The notation $|\cdot|$ shall be used to denote the size of an object,
whatever its nature. If $A$ is an measurable plane set then $|A|$ is its 
Lebesgue measure, also called its area ;
if $v$ is a vector $|v|$ is its $L^1$ norm ; if $\gamma$ is a rectifiable curve,
$|\gamma|$ is its $L^1$ length. We denote the boundary of a set using $\partial$.

By a \emph{domain} of the plane, we mean the closure
of the bounded component of a Jordan curve. In particular,
domains are compact and connected.
All rectangles and squares considered are assumed to have their sides parallel
to the coordinate axes. The square centered at $0$
with side length $2\lambda$ is denoted by $B_\infty(\lambda)$: it is the $\lambda$
ball of the $L^\infty$ metric. Squares are known to minimize
$L^1$ perimeter among plane domains of given area. 

The measure of the distance between compact plane sets $A,B$ we use in
our main result is the $L^\infty$ Haussdorf metric :
$$d_\infty(A,B)=\inf\{\lambda\leqslant0\,|\, A\subset B+B_\infty(\lambda)\mbox{ and }
  B\subset A+B_\infty(\lambda)\}.$$
Let us explain why this metric is natural here. One way to prove that
almost isoperimetric domains are close to minimizers is to prove
that they contain a minimizer of radius $r$ and are included in another of radius
$R$, with small radii difference $R-r$ and same center.
In the euclidean space, such inclusions
imply that the considered domain is at Haussdorf distance at most $(R-r)/2$
from some ball. However, balls and minimizers are different in the $L^1$ plane, 
so that if $A$ is between concentric squares of radii $R$ and $r$, one can only say that
it is at $L^1$ Haussdorf distance $R-r$ from some square, while the
$L^\infty$ Hausdorf distance bound is the expected $(R-r)/2$.

It would certainly be possible to use the $L^1$ Haussdorf metric,
and we expect that arguments of the same kind that those we use to
prove Theorem \ref{theo:main}, but more involved, would give a constant better
than the $1/16$ obtained using the inequality $d_1\leqslant 2d_\infty$ and Theorem \ref{theo:main}.

\begin{theo}\label{theo:main}
Let $A$ be a domain of the $L^1$ Minkowski plane whose boundary is a rectifiable
curve, and assume that
\begin{equation}
|\partial A|^2\leqslant (16+\varepsilon)|A|.\label{eq:hyp}
\end{equation}
Then there is a square $S$ such that
\begin{equation}
d_\infty(A,S)^2\leqslant\frac{\varepsilon |A|}{64}. \label{eq:conc}
\end{equation}
\end{theo}
We shall also see that Theorem \ref{theo:main} is sharp and show that up to $L^1$
isometry and homothety the domains that achieve the bound are the rectangles
and the squares with one square deleted at a corner.

A second possible measure of the distance between domains of the same area,
which present the advantage to be suitable to higher dimension as well,
is simply the gap between their area and that of their intersection.  
In this respect we prove the following.
\begin{prop}\label{prop:area}
Let $A$ be a domain of the $L^1$ Minkowski plane whose boundary is a rectifiable
curve, and assume that \eqref{eq:hyp} holds with $\varepsilon$ sufficiently small.
Then there is a square $S$ such that $|S|=|A|$ and:
\begin{equation}
|S\cap A| \geqslant (1-\frac{\sqrt{\varepsilon}}{4}+O(\varepsilon))|A|.\label{eq:conc2}
\end{equation}
In terms of Fraenkel asymmetry, this reads:
\begin{equation}
\frac{|S\Delta A|}{|A|} \leqslant \frac{\sqrt{\varepsilon}}2+O(\varepsilon)
\label{eq:asymmetry}
\end{equation}
\end{prop}
We shall see that the $1/2$ constant in \eqref{eq:asymmetry}
is sharp.

Surprisingly enough, it seems that these results are new, although 
similar ones can be deduced from the much more general
\cite{Figalli-Maggi-Pratelli} (but with a non-optimal constant)
and \cite{Peri-Wills-Zucco} (only when $A$ is convex).

\section{Proof of the inequalities}

Assume $A$ satisfies \eqref{eq:hyp} for some $\varepsilon$
and let $R$ be the smallest rectangle containing $A$.
This rectangle plays the role of a convex hull.

\begin{lemm}\label{lemm:1}
We have $|\partial A|\geqslant |\partial R|$.
\end{lemm}

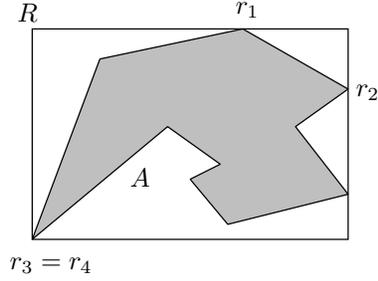
\begin{figure}[htp]\begin{center}
\input{rect.pstex_t}
\caption{The $L^1$ perimeter of $A$ is at least that of $R$}\label{fig:rect}
\end{center}\end{figure}

\begin{proof}
Since $R$ is minimal, each of its sides contains a point of the
boundary of $A$. Denote $r_1,r_2,r_3,r_4$ such points so that
$r_i$ and $r_{i+1}$ lie on two adjacent sides of $R$ for all
$i$ (modulo $4$). It is possible that some $r_i=r_{i+1}$, but this
does not affect what follows.

There are four curves $\gamma_i$ in $\partial A$ that connect $r_i$
to $r_{i+1}$ and meet only at their endpoints (see figure \ref{fig:rect}).
Similarly, the boundary of $R$ is made of four curves $\eta_i$ connecting
$r_i$ to $r_{i+1}$. Since $R$ is a rectangle, the $\eta_i$ are $L^1$ geodesics.
The length of $\gamma_i$ is at least $|r_i-r_{i+1}|=|\eta_i|$, so that 
$$|\partial A|=|\gamma_1|+|\gamma_2|+|\gamma_3|+|\gamma_4|\geqslant 
|\eta_1|+|\eta_2|+|\eta_3|+|\eta_4|=|\partial R|.$$
\end{proof}

Let $\ell$ and $\alpha$ be such that 
$\ell-2\alpha$ and $\ell+2\alpha$ are the side lengths of $R$.

\begin{lemm}\label{lemm:alpha} 
We have
\begin{equation}
|A|\leqslant \ell^2 \leqslant \frac{16+\varepsilon}{16}|A| \label{eq:ell}
\end{equation}
and
\begin{equation}
\alpha^2 \leqslant \frac{\varepsilon|A|}{64} \label{eq:alpha}
\end{equation}
\end{lemm}

\begin{proof}
From previous lemma we have $|\partial A|\geqslant 4\ell$, so that using 
\eqref{eq:hyp} we get $16\ell^2 \leqslant (16+\varepsilon)|A|$.
Since $A\subset R$ we have $|A|\leqslant|R|= \ell^2-4\alpha^2$
and \eqref{eq:ell} follows.

Next we have
\begin{eqnarray*}
16\ell^2 &\leqslant& (16+\varepsilon)(\ell^2-4\alpha^2) \\
0        &\leqslant& \varepsilon\ell^2-4(16+\varepsilon)\alpha^2 \\
\alpha^2 &\leqslant& \frac{\varepsilon\ell^2}{4(16+\varepsilon)}\\
\alpha^2 &\leqslant& \frac{\varepsilon}{64}|A|
\end{eqnarray*}
and we are done.
\end{proof}

Note that this lemma
is sufficient to deduce the $L^1$ isoperimetric inequality and its equality
case: if $\varepsilon=0$, then $\alpha=0$ and $|A|=\ell^2$.

\subsection{Proof of Theorem \ref{theo:main}}

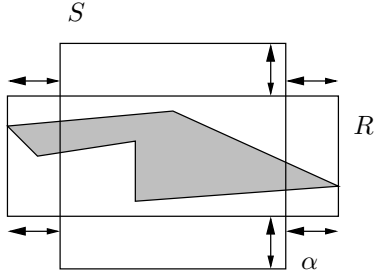
\begin{figure}[htp]\begin{center}
\input{central_square.pstex_t}
\caption{The closest square to $R$.}\label{fig:central_square}
\end{center}\end{figure}

We have $\min_S d_\infty(A,S)\geqslant\alpha$ (where $S$
runs over all squares, see figure \ref{fig:central_square})
and if there is equality, Lemma \ref{lemm:alpha} is sufficient
to conclude.
We therefore assume $\delta:=\min_S d_\infty(A,S)>\alpha$.

The following is the main step of the proof.
\begin{lemm}\label{lemm:2}
We have either
$$|A|\leqslant \ell^2-4\alpha^2-8\delta(\delta-\alpha)$$
or
$$|A|\leqslant \ell^2-4\alpha^2-4\delta^2.$$
\end{lemm}

\begin{proof}
Choose the origin so that $R$ has its bottom side at height $0$.
Let $S_\eta$ be the square that is at distance $\delta$
from each short side of $R$ (so that it has side length $\ell+2\alpha-2\delta$)
and whose bottom side is at height $\eta$.

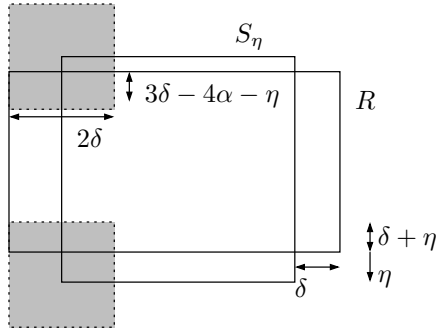
\begin{figure}[htp]\begin{center}
\input{square2.pstex_t}
\caption{The domain $A$ avoids one of the grey squares.}\label{fig:square2}
\end{center}\end{figure}

When $\eta\in[-\delta,3\delta-4\alpha]$, $R$ and $A$ are contained in the
$L^\infty$ neighborhood of size $\delta$ around $S_\eta$, thus
there is some point $p_\eta\in S_\eta$ that is at $L^\infty$ distance
at least $\delta$ from $A$.

This excludes $A$ from a square centered at $p_\eta$; the worst case
(with respect to our goal of bounding $|A|$ from above) is when this
excluded squares intersect only small parts of $R$ and have maximal
overlap. This is achieved when $p_\eta$ is a corner of $S_\eta$ for all $\eta$
and the short side of $R$ closest to $p_\eta$ is constant.

In this case, for each $\eta$, if $p_\eta$ is a lower corner
then there is a $2\delta\times(\delta+\eta)$ sub-rectangle of $R$ excluded,
else $p_\eta$ is a upper corner and there is a $2\delta\times(3\delta-4\alpha-\eta)$
sub-rectangle of $R$ excluded. These values assume that $\eta<\delta$ and
$3\delta-4\alpha-\eta<\delta$ respectively, otherwise there is simply an
excluded square of area $4\delta^2$.

Let $x$ be the supremum of the $\eta$ such that $p_\eta$ is a lower corner.
There is an excluded sub-rectangle of area $2\delta\times(\delta+x)$
and for all $\eta>x$ the point $p_\eta$ must be a higher corner,
so that there is another excluded sub-rectangle of area $2\delta\times(3\delta-4\alpha-x)$.

Summing up, either there are excluded sub-rectangles of total area at least
$2\delta\times4(\delta-\alpha)$, or there is an excluded sub-rectangle 
of area at least $4\delta^2$, and we get the desired bounds on
$|A|$.
\end{proof}

We can now conclude the proof of Theorem \ref{theo:main}.
First, if
$|A|\leqslant (\ell^2-4\alpha^2)-8\delta(\delta-\alpha)$
then we have
\begin{eqnarray*}
|\partial A|^2 &\leqslant& (16+\varepsilon)|A| \\
(4\ell)^2               &\leqslant& (16+\varepsilon)\big(\ell^2-4\alpha^2
	         -8\delta(\delta-\alpha)\big)\\
0              &\leqslant& \varepsilon \ell^2-4(16+\varepsilon)
                 (\alpha^2+2\delta(\delta-\alpha))\\
\alpha^2+2\delta(\delta-\alpha) &\leqslant& \frac{\varepsilon\ell^2}{4(16+\varepsilon)}\\
2\delta^2-2\alpha\delta+\alpha^2 &\leqslant& \frac{\varepsilon}{64}|A|
\end{eqnarray*}
the last inequality coming from \eqref{eq:ell}.

Since the function $x\mapsto2\delta^2-2x\delta+x^2$ is minimal when $x=\delta$,
we have 
$2\delta^2-2\alpha\delta+\alpha^2\geqslant 2\delta^2-2\delta^2+\delta^2=\delta^2$,
so that $\delta^2\leqslant \frac{\varepsilon}{64}|A|$.

In the case when $|A|\leqslant \ell^2-4\alpha^2-4\delta^2$, we get:
\begin{eqnarray*}
|\partial A|^2 &\leqslant& (16+\varepsilon)(\ell^2-4\alpha^2-4\delta^2) \\
16\ell^2 &\leqslant& 16\ell^2+\varepsilon\ell^2-4(16+\varepsilon)(\alpha^2+\delta^2)\\
\alpha^2+\delta^2 &\leqslant& \frac{\varepsilon\ell^2}{4(16+\varepsilon)}\\
\delta^2 &\leqslant&\frac{\varepsilon}{64}|A|
\end{eqnarray*}

\subsection{Proof of Proposition \ref{prop:area}}

Let $\mu=\max_S|S\cap A|/|A|$ where $S$ runs over the squares
having same area than $A$.

\begin{lemm}
One of the following holds:
\begin{eqnarray*}
\mu &\geqslant& 2-\frac{\ell^2-4\alpha^2}{|A|}\\
\mu &\geqslant& 2-\frac{\ell+2\alpha}{\sqrt{|A|}}.
\end{eqnarray*}
\end{lemm}

\begin{figure}[htp]\begin{center}
\input{squarebis.pstex_t}
\caption{The domain $A$ is included in $R$ and cannot meet a too large
   proportion of $S_0$}\label{fig:squarebis}
\end{center}\end{figure}
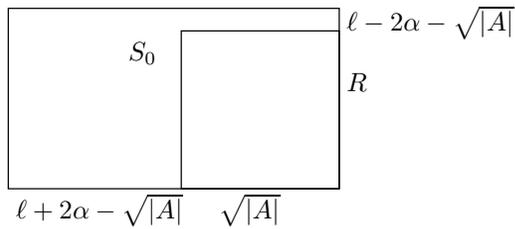

\begin{proof}
Define $S_0$ to be a square that shares a corner of $R$ and
intersects its interior, and that have the same area than $A$
(see figure \ref{fig:squarebis}).
The definition of $\mu$
implies that $|A\cap S_0|\leqslant \mu|A|$.

If $\sqrt{|A|}\leqslant \ell-2\alpha$ we have:
\begin{eqnarray*}
|A| &\leqslant& \mu|A|+(\ell+2\alpha)(\ell-2\alpha-\sqrt{|A|})
                +\sqrt{|A|}(\ell+2\alpha-\sqrt{|A|})\\
    &\leqslant& \mu |A|+\ell^2-4\alpha^2-\sqrt{|A|}(\ell+2\alpha)
                +\sqrt{|A|}(\ell+2\alpha)-|A|\\
    &\leqslant& \ell^2-4\alpha^2+(\mu-1)|A|\\
\mu|A| &\geqslant& 2|A|-(\ell^2-4\alpha^2).
\end{eqnarray*}
Otherwise, we get
\begin{eqnarray*}
|A| &\leqslant& \mu|A|+(\ell-2\alpha)(\ell+2\alpha-\sqrt{|A|}) \\
    &\leqslant& \mu|A|+\sqrt{|A|}(\ell+2\alpha-\sqrt{|A|})\\
\mu|A| &\geqslant& 2|A|-(\ell+2\alpha)\sqrt{|A|}
\end{eqnarray*}
\end{proof}

If the first conclusion holds, using Lemma \ref{lemm:alpha} it comes
\begin{eqnarray*}
\mu &\geqslant& 2-\frac{(\ell^2-4\alpha^2)(16+\varepsilon)}{16\ell^2}\\
         &\geqslant&2-\frac{16+\varepsilon}{16}=1-\frac{\varepsilon}{16}
\end{eqnarray*}

If the second conclusion holds, using Lemma \ref{lemm:alpha} we get
\begin{eqnarray*}
\mu &\geqslant& 2-\frac{\ell}{\sqrt{|A|}}-\frac{2\alpha}{\sqrt{|A|}}\\
    &\geqslant& 2-\sqrt{1+\frac{\varepsilon}{16}}-\frac{\sqrt{\varepsilon}}4\\
    &\geqslant& 1-\frac{\sqrt{\varepsilon}}4+O(\varepsilon)
\end{eqnarray*}
But for all sufficiently small $\varepsilon$, this second expression is smaller than
$1-\varepsilon/16$, and Proposition \ref{prop:area} is proved.

\section{Sharpness}\label{sec:sharpness}

Two examples showing sharpness of Theorem \ref{theo:main} steam
out from its proof. 

\begin{figure}[htp]\begin{center}
\includegraphics{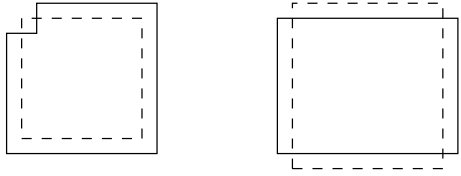}
\caption{Two domain that are almost isopermetric and as far as possible
  from squares: a square with a small corner deleted and a rectangle with 
  sides of almost the same lengths.}
\end{center}\end{figure}

The first one is the domain $S'_\delta$ obtained from the unit square
by deleting a $2\delta\times 2\delta$ square at one corner ($\delta<1/2)$).
We have $|S'_\delta|=1-4\delta^2$ and $|\partial S'_\delta|=4$,
so that \eqref{eq:hyp} holds with 
$$\varepsilon=\frac{64\delta^2}{1-4\delta^2}$$
and $\inf_S d_\infty(S,S'_\delta)=\delta$ so that equality holds in \eqref{eq:conc}.

The second one is
the rectangle $R_\alpha$ whose side length
are $1-2\alpha$ and $1+2\alpha$ (where $\alpha<1/2)$).
We have $|R_\alpha|=1-4\alpha^2$, $|\partial R_\alpha|=4$
and $\inf_S d_\infty(S,R_\alpha)=\alpha$ so that
\eqref{eq:conc} is an equality once again.

Let us show that $S'_\delta$ and $R_\alpha$ are the only possible (up to homothety
and $L^1$ isometry) exemple realizing equality in both \eqref{eq:hyp}
and \eqref{eq:conc} for the same $\varepsilon$. In the first
case of Lemma \ref{lemm:2}, for $2\delta^2-2\alpha\delta+\alpha^2\leqslant\delta^2$
to be an equality it is necessary that $\alpha=\delta$, so that $A$ must be equal to $R$ (otherwise
$R$ would have smaller isoperimetric inequality and same distance to squares).
In the second case of the lemma, one is lead to $\alpha=0$ in the last lines
of the proof of Theorem \ref{theo:main}, so that $R$ is a square and according
to the proof of Lemma \ref{lemm:2}, $A$ is contained in a $S'_\delta$ having
the same isoperimetric deficit and the same minimal rectangle. They must therefore
be equal.

At last, $R_\alpha$ shows asymptotic sharpness of Proposition \ref{prop:area}:
$$\sup_{|S|=|R_\alpha|} |S\cap R_\alpha|=(1-2\alpha)\sqrt{1-4\alpha^2}=1-2\alpha+o(\alpha)$$
and 
$$1-\frac14\sqrt{\varepsilon}=1-2\alpha+o(\alpha)$$
when $\varepsilon$ takes the extremal value $64\alpha^2/(1-4\alpha^2)$.

\bibliographystyle{alpha}
\bibliography{biblio}

\end{document}

%% file: rect.pstex_t
\begin{picture}(0,0)%
\includegraphics{rect.pstex}%
\end{picture}%
\setlength{\unitlength}{4144sp}%
\begingroup\makeatletter\ifx\SetFigFontNFSS\undefined%
\gdef\SetFigFontNFSS#1#2#3#4#5{%
  \reset@font\fontsize{#1}{#2pt}%
  \fontfamily{#3}\fontseries{#4}\fontshape{#5}%
  \selectfont}%
\fi\endgroup%
\begin{picture}(2100,1741)(3901,-3995)
\put(5266,-2401){\makebox(0,0)[lb]{\smash{{\SetFigFontNFSS{10}{12.0}{\rmdefault}{\mddefault}{\updefault}{\color[rgb]{0,0,0}$r_1$}%
}}}}
\put(3916,-3931){\makebox(0,0)[lb]{\smash{{\SetFigFontNFSS{10}{12.0}{\rmdefault}{\mddefault}{\updefault}{\color[rgb]{0,0,0}$r_3=r_4$}%
}}}}
\put(5986,-2896){\makebox(0,0)[lb]{\smash{{\SetFigFontNFSS{10}{12.0}{\rmdefault}{\mddefault}{\updefault}{\color[rgb]{0,0,0}$r_2$}%
}}}}
\put(3961,-2446){\makebox(0,0)[lb]{\smash{{\SetFigFontNFSS{10}{12.0}{\rmdefault}{\mddefault}{\updefault}{\color[rgb]{0,0,0}$R$}%
}}}}
\put(4636,-3436){\makebox(0,0)[lb]{\smash{{\SetFigFontNFSS{10}{12.0}{\rmdefault}{\mddefault}{\updefault}{\color[rgb]{0,0,0}$A$}%
}}}}
\end{picture}%

%% file: central_square.pstex_t
\begin{picture}(0,0)%
\includegraphics{central_square.pstex}%
\end{picture}%
\setlength{\unitlength}{4144sp}%
\begingroup\makeatletter\ifx\SetFigFontNFSS\undefined%
\gdef\SetFigFontNFSS#1#2#3#4#5{%
  \reset@font\fontsize{#1}{#2pt}%
  \fontfamily{#3}\fontseries{#4}\fontshape{#5}%
  \selectfont}%
\fi\endgroup%
\begin{picture}(2097,1672)(4039,-3950)
\put(6121,-3076){\makebox(0,0)[lb]{\smash{{\SetFigFontNFSS{10}{12.0}{\rmdefault}{\mddefault}{\updefault}{\color[rgb]{0,0,0}$R$}%
}}}}
\put(5806,-3886){\makebox(0,0)[lb]{\smash{{\SetFigFontNFSS{10}{12.0}{\rmdefault}{\mddefault}{\updefault}{\color[rgb]{0,0,0}$\alpha$}%
}}}}
\put(4411,-2401){\makebox(0,0)[lb]{\smash{{\SetFigFontNFSS{10}{12.0}{\rmdefault}{\mddefault}{\updefault}{\color[rgb]{0,0,0}$S$}%
}}}}
\end{picture}%

%% file: square2.pstex_t
\begin{picture}(0,0)%
\includegraphics{square2.pstex}%
\end{picture}%
\setlength{\unitlength}{4144sp}%
\begingroup\makeatletter\ifx\SetFigFontNFSS\undefined%
\gdef\SetFigFontNFSS#1#2#3#4#5{%
  \reset@font\fontsize{#1}{#2pt}%
  \fontfamily{#3}\fontseries{#4}\fontshape{#5}%
  \selectfont}%
\fi\endgroup%
\begin{picture}(2232,1959)(4039,-4393)
\put(6121,-3076){\makebox(0,0)[lb]{\smash{{\SetFigFontNFSS{10}{12.0}{\rmdefault}{\mddefault}{\updefault}{\color[rgb]{0,0,0}$R$}%
}}}}
\put(5401,-2671){\makebox(0,0)[lb]{\smash{{\SetFigFontNFSS{10}{12.0}{\rmdefault}{\mddefault}{\updefault}{\color[rgb]{0,0,0}$S_\eta$}%
}}}}
\put(4456,-3301){\makebox(0,0)[lb]{\smash{{\SetFigFontNFSS{10}{12.0}{\rmdefault}{\mddefault}{\updefault}{\color[rgb]{0,0,0}$2\delta$}%
}}}}
\put(6256,-4111){\makebox(0,0)[lb]{\smash{{\SetFigFontNFSS{10}{12.0}{\rmdefault}{\mddefault}{\updefault}{\color[rgb]{0,0,0}$\eta$}%
}}}}
\put(5761,-4201){\makebox(0,0)[lb]{\smash{{\SetFigFontNFSS{10}{12.0}{\rmdefault}{\mddefault}{\updefault}{\color[rgb]{0,0,0}$\delta$}%
}}}}
\put(6256,-3886){\makebox(0,0)[lb]{\smash{{\SetFigFontNFSS{10}{12.0}{\rmdefault}{\mddefault}{\updefault}{\color[rgb]{0,0,0}$\delta+\eta$}%
}}}}
\put(4861,-3031){\makebox(0,0)[lb]{\smash{{\SetFigFontNFSS{10}{12.0}{\rmdefault}{\mddefault}{\updefault}{\color[rgb]{0,0,0}$3\delta-4\alpha-\eta$}%
}}}}
\end{picture}%

%% file: squarebis.pstex_t
\begin{picture}(0,0)%
\includegraphics{squarebis.pstex}%
\end{picture}%
\setlength{\unitlength}{4144sp}%
\begingroup\makeatletter\ifx\SetFigFontNFSS\undefined%
\gdef\SetFigFontNFSS#1#2#3#4#5{%
  \reset@font\fontsize{#1}{#2pt}%
  \fontfamily{#3}\fontseries{#4}\fontshape{#5}%
  \selectfont}%
\fi\endgroup%
\begin{picture}(2052,1336)(4039,-4175)
\put(6076,-3346){\makebox(0,0)[lb]{\smash{{\SetFigFontNFSS{10}{12.0}{\rmdefault}{\mddefault}{\updefault}{\color[rgb]{0,0,0}$R$}%
}}}}
\put(6076,-2986){\makebox(0,0)[lb]{\smash{{\SetFigFontNFSS{10}{12.0}{\rmdefault}{\mddefault}{\updefault}{\color[rgb]{0,0,0}$\ell-2\alpha-\sqrt{|A|}$}%
}}}}
\put(5311,-4111){\makebox(0,0)[lb]{\smash{{\SetFigFontNFSS{10}{12.0}{\rmdefault}{\mddefault}{\updefault}{\color[rgb]{0,0,0}$\sqrt{|A|}$}%
}}}}
\put(4096,-4111){\makebox(0,0)[lb]{\smash{{\SetFigFontNFSS{10}{12.0}{\rmdefault}{\mddefault}{\updefault}{\color[rgb]{0,0,0}$\ell+2\alpha-\sqrt{|A|}$}%
}}}}
\put(4771,-3166){\makebox(0,0)[lb]{\smash{{\SetFigFontNFSS{10}{12.0}{\rmdefault}{\mddefault}{\updefault}{\color[rgb]{0,0,0}$S_0$}%
}}}}
\end{picture}%